# Development of the Measure of Assessment Self-Efficacy (MASE) for Quizzes and Exams

Kaitlin Riegel[a]*, Tanya Evans[b], and Jason M. Stephens[c]

[a,b]Department of Mathematics, University of Auckland, Auckland, NZ;

[c]Faculty of Education and Social Work, University of Auckland, Auckland, New Zealand;

*[a]krie235@aucklanduni.ac.nz (ORCID: 0000-0002-8187-2016);

[b]t.evans@auckland.ac.nz (ORCID: 0000-0001-5126-432X);

[c]jm.stephens@auckland.ac.nz (ORCID: 0000-0003-2404-6321)

# Development of the Measure of Assessment Self-Efficacy (MASE) for Quizzes and Exams

Kaitlin Riegel[a]*, Tanya Evans[b], and Jason M. Stephens[c]

[a,b]Department of Mathematics, University of Auckland, Auckland, NZ;

[c]Faculty of Education and Social Work, University of Auckland, Auckland, New Zealand;

Self-efficacy is a significant construct in education due to its predictive relationship with achievement. Existing measures of assessment-related self-efficacy concentrate on students' beliefs about content-specific tasks but omit beliefs around assessment-taking. This research aimed to develop and test the Measure of Assessment Self-Efficacy (MASE), designed to assess two types of efficacy beliefs related to assessment (i.e., 'comprehension and execution' and 'emotional regulation') in two scenarios (i.e., a low-stakes online quiz and a high-stakes final exam). Results from confirmatory factor analysis in Study 1 ($N = 301$) supported the hypothesised two-factor measurement models for both assessment scenarios. In Study 2, results from MGCFA ($N = 277$) confirmed these models were invariant over time and provided evidence for the scales' validity. Study 3 demonstrated the exam-related MASE was invariant across cohorts of students ($N$s = 277; 329). Potential uses of the developed scales in educational research are discussed.

**Introduction**

Bandura (1997) defined self-efficacy as "beliefs in one's capabilities to organize and execute the courses of action required to produce given attainments" (p. 3) In its many operationalisations across numerous studies, self-efficacy has been shown to be a significant predictor of academic achievement (Bandura, 2010; Ferla et al., 2009; Pajares & Graham, 1999; Richardson et al., 2012; Zimmerman, 2000). In educational contexts, experiences of success positively influence the development of self-efficacy while experiences of failure impair it (Usher & Pajares, 2009). In particular, this suggests success or failure in one assessment may influence a student's efficacy beliefs about future assessments. With this in mind, it is important to have an instrument to capture the dynamic efficacy beliefs with regards to taking assessments. This paper presents three sequential studies, which sought to develop and validate a measure designed to assess students' assessment-related efficacy beliefs.

Such a measure could be particularly useful in STEM fields (science, technology, engineering, and mathematics), and for mathematics in particular, which have notoriously high drop-out rates associated with poor performance (Chen, 2013). This is compounded by the growing concern about the students who opt out of pursuing mathematics at a tertiary level before they finish secondary school (Zacharopoulos et al., 2021). Therefore, it is important for researchers to have valid and reliable instruments to measure the relevant constructs in understanding causes of poor performance. Moreover, as research expands to include other constructs related to self-efficacy such as achievement emotions, motivation, and attitudes, it is critical to have appropriate measures, given the considerable impact such variables have on students' participation in mathematics and, consequently, tertiary STEM education and when there are equity issues at stake (Chen, 2013; Ma & Willms, 1999; Zacharopoulos et al., 2021). Furthermore, as mathematics education researchers develop new assessment innovations and interventions, the need to be able to accurately evaluate their usefulness from psychometric perspective becomes increasingly important.

*Assessment self-efficacy*

To date, many studies have only looked at general academic or domain-specific self-efficacy without including measures for assessment-related self-efficacy. Yet, it is common to hear students lament that they are "not test people," or "I know I can do it, just not on a test". Such comments suggest that measuring a student's belief in their abilities in mathematics may not be representative of their belief in their ability to achieve on assessment and a context-specific measure is needed for research aiming to address questions around, for example, performance, learning, attrition, or affect.

Bandura (1997) proposed that self-efficacy measures must be both context and content specific. In fact, subsequent research has demonstrated that context-related self-efficacy and content-related self-efficacy are separate constructs. Marsh et al. (2019) found generalised mathematics self-efficacy was distinct from test-related self-efficacy. Nielsen et al. (2017) argue that the Motivated Strategies for Learning Questionnaire (Pintrich et al., 1991) self-efficacy subscale (MSLQ-SE) consists of two context-specific scales: specific academic learning self-efficacy (SAL-SE) and specific academic exam self-efficacy (SAE-SE). According to Nielsen and Moore (2003):

> This contextual differentiation has been largely ignored in the literature, but it may be important because of the variability among students…These differences in students' reactions to, and performances in, different testing conditions may be related to their perceptions of their specific and context task efficacy. (p. 130)

Nielsen and Moore go on to report the significant difference they perceived in responses to their class Mathematics Self-Efficacy Scale (MSES) and test MSES, underscoring the need for context specificity when researching self-efficacy. The evidence suggests that context matters when investigating self-efficacy, and therefore research into self-efficacy as it pertains to assessment should also specifically consider students beliefs about their abilities in taking such an assessment.

A further barrier to accurate measurement of assessment-related self-efficacy is that the existing scales often measure outcome expectancies rather than assessment-related self-efficacy, which Bandura (2006) identifies as a separate construct. For example, all of the items on the SAE-SE (Nielsen et al., 2017) are outcome expectancies (e.g., I expect to do well in this class). Similarly, the Mathematics self-efficacy scale of Pietsch, Walker, and Chapman (2003) consists of items asking about achieving certain grades (e.g., I am able to achieve at least 50% in my mathematics course this year). An exception is Coohey and Cummings (2019), who developed the Confidence in Test-Taking Skills (CTTS), aimed at assessing student self-efficacy in multiple-choice exams. The five statements were based on student beliefs about their understanding and abilities on a multiple-choice exam, such as their ability to make an educated guess. However, this measure caters to the specific assessment so cannot be used for a different format.

Despite how assessment-related self-efficacy is currently investigated, the beliefs students hold around taking an assessment are not limited to its content. There are other features, such as beliefs about their ability to manage their time and regulate their emotions in a specific environment. So far, most research has investigated task-based self-efficacy that was tied to a certain context, with items such as "I can calculate percentages in class/in a test", but little has sought to develop measures based on the broad range of aspects involved in taking an assessment, which students may or may not be efficacious about.

*The present investigation*

The purpose of the present study was to develop and validate the Measure of Assessment Self-Efficacy (MASE). To do so, we conducted three (sequential) studies. Study 1 aimed to develop items for a scale under two assessment scenarios and validate the latent efficacy factors. Incorporating beliefs around both cognitive and affective aspects of assessment-taking into the item development presented an opportunity for conception of a robust measure of assessment-related self-efficacy that accounts for multi-faceted essence of the construct. Study 2 sought to assess the invariance of the scale for two assessment scenarios across time and to provide evidence of construct validity through employing the Fornell-Larcker criterion (Fornell & Larcker, 1981), as well as examining correlations with achievement emotions and academic performance. Previous research has demonstrated that self-efficacy is associated positively with positive emotions and negatively with negative emotions (Pekrun & Stephens, 2010; Pekrun et al., 2011; Luo et al., 2016). In particular, Pekrun et al. (2004) argued that self-efficacy relates to positive test-related emotions. It has also been shown that high anxiety can undermine self-efficacy (Usher & Pajares, 2009). Therefore, the convergent and divergent validity of the scales was investigated with these relationships in mind and the known association between self-efficacy and academic achievement. Finally, Study 3 examined the invariance of the measure for two assessment scenarios between cohorts of students.

**Study 1**

*Purpose*

The aim of this study was to develop items for a two-factor model for assessment-related self-efficacy and test the model in two distinct assessment scenarios.

*Participants and procedure*

Study 1 used a cross-sectional research design involving the completion of an anonymous survey by a sample of undergraduate university students ($N = 301$) procured through Prolific (a crowdsourcing platform based in the UK that, like MTurk, pays individuals a nominal fee for participating in research projects). The sample consisted of 193 students who identified as male, 104 as female, three as gender diverse, and one who declined to answer. Out of the respondents, 62% identified as Caucasian, 19% as Hispanic or Latino, 4% as Asian, 4% as Middle Eastern, 1% as Black or African-American, 7% as from another ethnicity, and 3% as being of multiple ethnicities. One participant did not respond to this question. An approval to conduct the study from the University of Auckland Human Participants Ethics Committee was granted (approval number 024726).

*Measures*

*The Measure of Assessment Self-Efficacy (MASE)*

10 original items were created to assess two latent factors associated with participants' beliefs in their ability to understand, perform, and regulate emotions before and during an assessment. The first factor, *comprehension and execution*, aimed to measure participants' beliefs around their abilities to understand and execute on an assessment while the second factor, *emotional regulation*, aimed to measure participants' beliefs around their capacity to regulate their emotions as they prepared for and completed an assessment. The items are presented and discussed in the Results section.

      Following Bandura's (2006) recommendations, the items were written to measure participants perceived capabilities and responses were collected using a slider scale from 0 to 100 (where 0 = *Cannot do at all*, 50 = *Moderately sure can do*, and 100 = *Highly certain can do*). Each participant responded to the self-efficacy items in both of two assessment

scenarios, a low-stakes online quiz and a high-stakes final exam, which were randomly presented.

> Mathematics Quiz Scenario
> Imagine that you've enrolled in a mathematics course that has a QUIZ worth 0.5% of your final grade. The quiz is open book and taken online at any time you choose in a 24-hour period. You have 30 minutes to answer two multiple-choice questions based on the previous lecture. The quiz allows for two attempts and your best result is recorded.
>
> Mathematics Exam Scenario
> Imagine that you've enrolled in a mathematics course that has a final EXAM worth 50% of your final grade. The exam contains short and long answer questions. The exam is invigilated and is two hours long.

These two assessment scenarios were chosen because of their intrinsic differences as assessment methods. Further, high-risk final exams are still a standard assessment approach in tertiary mathematics, while forms of online assessment have been increasingly implemented, particularly with the occurrence of the COVID-19 pandemic. As such, the choice of these contrasting scenarios to validate scales for use in further research was deemed as reasonable.

### *Data analyses*

There was no missing data in this analysis. Normed chi-square values and several alternative fit indices were used to determine model fit: $\chi^2/df < 3.0$, TLI > .95, CFI > .95, RMSEA < .06, and SRMR < .05 for a "good" fit, and $\chi^2/df < 5.0$, TLI > .90, CFI > .90, RMSEA < .08, and SRMR < .08 for an "acceptable" fit (Hu & Bentler, 1999; Ullman & Bentler, 2003), though SRMR < .10 has also been considered as "acceptable" (MacCallum et al., 1996; Schermelleh-Engel et al., 2003). Factor loadings over .70 were considered good (Dash & Paul, 2021; Hair et al., 2022). Cronbach's alpha of over .70 were considered good (Taber, 2017). Analyses were conducted using version 27 of SPSS and its AMOS programme. Where applicable, the same standards and software were employed for all studies presented in this paper.

### *Results*

#### *Confirmatory factor analysis*

Confirmatory factor analysis demonstrated the ten-item, two-factor model provided an acceptable fit for assessment self-efficacy in scenarios: quiz ($\chi^2/df = 2.75$; TLI = .967; CFI = .975; RMSEA = .076, 90%CI[.058-.095]; SRMR = .0326) and exam ($\chi^2/df = 3.43$; TLI = .948; CFI = .961; RMSEA = .090, 90%CI[.072-.108]; SRMR = .0393). Appendix A outlines the distribution and correlations of the items in the scale for both scenarios. Specifically, with the exception of the slightly high RMSEA values, all indices suggested a *good* fit of the models to the data. An examination of the standardised residuals and modification indices suggested that item 8 could be a source of misspecification for both models. We further noted that modification indices suggested that item 6 in the quiz scenario and items 3 and 9 in the exam scenario could also be sources of misspecification for their respective models.

Table 1 provides a summary of the factors, items, and loadings for both the quiz and the exam. With the exception of item 3, all factor loadings were good to excellent (range = .73 to .89). Similarly, the Cronbach's alpha values for all four factors were also good to excellent (range = .86 to .91).

Table 1. Factors, items, and loadings for quiz and exam-related self-efficacy in Study 1.

| Factor/Item | Loading Quiz | Loading Exam |
|---|---|---|
| **Comprehension and execution** | | |
| 1. I can understand the content and skills needed for the assessment. | .87 | .80 |
| 2. I can fully master the requirements I need for this assessment. | .88 | .85 |
| 3. When preparing for the assessment, I can organize my time well. | .51 | .51 |
| 4. I can understand the questions in the assessment. | .85 | .81 |
| 5. I can correctly answer the questions on the assessment. | .89 | .86 |
| 6. During the assessment, I can answer the questions within the time constraints. | .76 | .73 |
| (Scale alpha) | (.90) | (.89) |
| **Emotional regulation** | | |
| 7. During my preparation, I am able to cope with my negative emotions toward the assessment. | .84 | .84 |
| 8. Even when I struggle while studying, I am able to stay positive about my ability to succeed. | .83 | .87 |
| 9. During the assessment, I am able to cope with my negative emotions toward the assessment. | .84 | .80 |
| 10. Even when I struggle during the assessment, I am able to stay positive about my ability to succeed. | .88 | .87 |
| (Scale alpha) | (.91) | (.86) |

*Note.* Loadings are standardised beta weights.

### Conclusions (Study 1)

Study 1 aimed to develop and validate the factor structure of the MASE under two assessment scenarios. Results from confirmatory factor analyses supported the two-factor models (assessing *comprehension and execution* and *emotional regulation*), indicating acceptable fits for both quiz-related and exam-related self-efficacy. The 10-item scale did not provide *good* fits due to the high RMSEA, however, the decision was made to note potential sources of misspecification for the models and proceed with the full set of items in order to prioritise identical and broad scales in both scenarios. The analysis also suggested that the model fit might be improved, particularly through the removal or rephrasing of item 8. Accordingly, as detailed below, this item was revised for use for Study 2.

## Study 2

### Purpose

The aim of this second study was to establish measurement invariance of the quiz and exam-related MASE over time using longitudinal data. We also sought to assess construct validity through use of the Fornell-Larcker criterion (Fornell & Larcker, 1981), as well as the correlations between the MASE factors, achievement emotions, and academic achievement.

### Participants

Study 2 involved a sample of university students enrolled in a standard second-year service mathematics course. Specifically, 410 students were enrolled in the course at the start of the semester, nine of whom later withdrew. Of the remaining 401 students, 379 provided consent to use their data from the course for research purposes. After data cleaning (e.g., straight-lining or missing a large number of responses, specifically, entire scales), time point one had 356 responses, time point two had 329, and time point three had 308. In total, there were 277 students who completed the survey at all three data collection points. Out of 242 items, ten participants missed one item, one participant missed six items, and one participant missed ten items. Missing data were inserted using EM-imputation. The dataset consisted of 143 participants who identified as male, 131 as female, two as gender diverse, and one who declined to respond. Out of the participants, 71% identified as Asian, 12% as NZ European, 2% as Other European, 1% as Māori or Pacific peoples, 1% as Middle Eastern, 3% as from another ethnicity, and 9% as being of multiple ethnicities. An approval to conduct the study from the University of Auckland Human Participants Ethics Committee was granted (approval number 024710).

*Procedure*

A questionnaire was distributed online in the first, seventh, and last (12th) week of the semester. Students were allocated a small amount (0.1%) of their final course grade for each questionnaire they completed. Importantly, the assessment structure for the course in this study included weekly online quizzes (which began during the second week of the semester), each worth only approximately 0.5% of the course grade, and a final exam (at the end of the semester) worth 50% of the course grade.

*Measures*

*The MASE*

The survey instrument included the MASE described in Study 1 with two exceptions. First, based on findings of Study 1, one of the 10 items (Item 8) was revised by deleting the introductory phrase "Even when I struggle…," as we suspected this item not fitting was a consequence of participants not widely feeling they struggle while studying. Second, in order to provide the students with more immediate context for the type of assessment they were being asked about, the first line of instructions was altered to read, "Imagine you are enrolled in a course like Maths 208". Each participant was again asked to respond to the scale under both of the previously described assessment scenarios.

*Achievement emotions*

To measure achievement emotions, the survey included an adapted version of the Academic Emotions Questionnaire (AEQ) (Pekrun et al., 2011). The test-related section of the AEQ uses five-point Likert-type response scale (1 = *Strongly Disagree* to 5 = *Strongly Agree*) to measure eight types of emotions related to academic achievement. For the purposes of the present investigation, the AEQ was adapted in several ways. First, we chose only to investigate emotions experienced *before* and *during* (and not *after*) an assessment to align with the temporal aspect of the MASE items. Secondly, it was shortened to assess only two positive and two negative emotions for the purposes of assessing convergent and divergent validity - enjoyment, hope, anxiety, and hopelessness, (eliminating anger, pride, and shame, which often feature more heavily *after* an assessment). Finally, because there were very few items measuring enjoyment and hope experienced *during* an exam, we developed four new items to strengthen the validity of these factors. Confirmatory factor analysis initially indicated an unacceptable fit for these models. Modification indices suggested several items as a source of misspecification, which were sequentially removed until acceptable fits were obtained. Results from our final CFAs indicated an acceptable fit for both models: before an exam ($\chi^2/df$ = 2.09, CFI = .92, RMSEA = .06, $n$ of items = 18) and during an exam ($\chi^2/df$ = 2.72, CFI = .90, RMSEA = .08, $n$ of items = 18). Please see Appendix B for all scales, their items and loadings, and the deleted items.

*Exam stress*

One original item was created to assess participants' level of stress related to a mathematics exam. Specifically, after reading the exam assessment scenario described above, participants were asked to use a nine-point Likert-type scale (1 = *Not stressful at all* to 9 = *Extremely stressful*) to respond to the question, "How stressful do you perceive this mathematics EXAM to be?"

    Finally, self-reported gender, prerequisite grades as a measure for prior achievement (1 = C- to 9 = A+), and final exam marks for the course were used in this analysis.

*Data analyses*

Multi-group confirmatory factor analysis (MGCFA) was conducted to test the measurement invariance of the MASE across its use at three time points in this study. Following Cheung and Rensvold (2002), the level of invariance achieved was determined by a change in CFI of < .01.

To establish convergent validity, we considered the average variance extracted (AVE) to be greater than .50 and composite reliability (CR) of above .70. To establish discriminant validity, we checked the √AVE for each of the factors was greater than the correlation involving the factors (Fornell & Larcker, 1981). Further, Pearson correlation coefficients were employed to assess convergent and divergent validity.

*Results*

*Confirmatory factor analysis*

Using the time point 1 data, we sought to confirm the ten-item MASE for both contexts, but models did not offer an acceptable fit for either the quiz ($\chi^2/df$ = 4.45; TLI = .924; CFI = .942; RMSEA = .112, 90%CI[.094-.130]; SRMR = .059) or the exam ($\chi^2/df$ = 5.70; TLI = .889; CFI = .916; RMSEA = .130, 90%CI[.113-.149]; SRMR =.057). Working through each model independently and using both modification indices and theoretical considerations, several successive alterations were made to improve model fit. Specifically, we allowed for the correlation of error terms for items seven and nine, as we perceive them theoretically to be the same statement situated before taking the assessment and during taking the assessment. We then removed items sequentially based on modification indices, including the removal of items identified as potential sources of misspecification in Study 1. Appendix C outlines the steps taken in improving model fit, with the intention of obtaining a good fit, for the quiz-related MASE (MASE-Q) and the exam-related MASE (MASE-E). The RMSEA for the MASE-Q was acceptable, though still did not indicate a good fit. However, we decided to proceed with this model as all other fit indices were good and we did not want to cut further items from the scale. The MASE-E demonstrated a good model fit across all fit indices.

As detailed in Table 2, after completing modifications to improve fit, the MASE-Q and MASE-E each had 7-items. For both scales, four items measured *comprehension and execution* and three measured *emotional regulation*. Items 1, 2, and 5 were shared on both scales for the first factor (*comprehension and execution*), while item 3 loaded only onto the MASE-Q and item 4 only onto the MASE-E. Items 7 and 10 were shared on both scales for the second factor (*emotional regulation*), while item 8 loaded only onto the MASE-E and item 9 only onto the MASE-Q. With the exception of item 3 on the MASE-Q ($\lambda$ = .68), all factor loadings were good or excellent (range = .76 to .89). Finally, the Cronbach's alpha values for all four factors were also good to excellent (.84 to .91).

Table 2. Factors, items, and loadings for the MASE-Q and MASE-E in Study 2.

| Factor/item | Loading Quiz | Loading Exam |
|---|---|---|
| **Comprehension and execution** | | |
| 1. I can understand the content and skills needed for the assessment. | .83 | .86 |
| 2. I can fully master the requirements I need for this assessment. | .81 | .88 |
| 3. When preparing for the assessment, I can organize my time well. | .68 | - |
| 4. I can understand the questions in the assessment. | - | .82 |
| 5. I can correctly answer the questions on the assessment. | .86 | .83 |
| 6. During the assessment, I can answer the questions within the time constraints. | - | - |
| (Scale alpha) | (.87) | (.91) |
| **Emotional regulation** | | |
| 7. During my preparation, I am able to cope with my negative emotions toward the assessment. | .86 | .73 |
| 8. While studying, I am able to stay positive about my ability to succeed. | - | .89 |
| 9. During the assessment, I am able to cope with my negative emotions toward the assessment. | .80 | - |
| 10. Even when I struggle during the assessment, I am able to stay positive about my ability to succeed. | .85 | .76 |
| (Scale alpha) | (.89) | (.84) |

*Note.* Loadings are standardised beta weights.

Table 3 outlines the results of confirmatory factor analysis of the two, two-factor models at each time point, employing fit indices described previously to assess model fit. With the exception of the slightly high RMSEA values ($\geq$ .06), which may be a consequence of our small sample size and low degrees of freedom (Kenny et al., 2015), all indices suggested a "good" fit of the data to both measurement models. Specifically, all $\chi^2/df$ were "good" (< 3.0), as were the TLI and CFI values ($\geq$ .95) and SRMR (< .05).

Table 3. Two-factor models of assessment-related self-efficacy in Study 2.

| Model | $\chi^2$ | *df* | $\chi^2/df$ | TLI | CFI | RMSEA | (90% CI) | SRMR |
|---|---|---|---|---|---|---|---|---|
| **MASE-Q** | | | | | | | | |
| Time point 1 | 28.54* | 12 | 2.38 | .976 | .986 | .071 | (.037-.104) | .023 |
| Time point 2 | 27.25* | 12 | 2.27 | .984 | .991 | .068 | (.034-.102) | .019 |
| Time point 3 | 36.029** | 12 | 3.00 | .977 | .987 | .085 | (.054-.118) | .023 |
| **MASE-E** | | | | | | | | |
| Time point 1 | 19.14 | 13 | 1.47 | .992 | .995 | .041 | (.000-.078) | .024 |
| Time point 2 | 17.90 | 13 | 1.38 | .996 | .997 | .037 | (.000-.075) | .012 |
| Time point 3 | 35.77** | 13 | 2.75 | .982 | .989 | .080 | (.049-.111) | .016 |

*Note. N* = 277. *$p$ < .05, **$p$ < .005.

*Measurement invariance over time*

In light of the CFA results, we proceeded with testing of the two models for measurement invariance across three time points using MGCFA. Table 4 outlines the levels of invariance achieved over time for both the MASE-Q and MASE-E. Based on the observed changes in

CFI (< .01), the MASE-Q demonstrated "strict" invariance and the MASE-E "strong" invariance. Importantly, only "strong" (i.e, configural, metric, and scalar) invariance is required to conclude model equivalence and thus allow for comparisons of latent factor means across time.

Table 4. Measurement invariance of the MASE-Q and MASE-E within a cohort.

| Model | CFI Unconstrained | RMSEA Unconstrained | ΔCFI (*ΔRMSEA*) | | |
|---|---|---|---|---|---|
| | | | Metric Invariance | Scalar Invariance | Residual Invariance |
| MASE-Q | .988 | .043 | .000 (*.005*) | .009 (*.009*) | .006 (*.000*) |
| MASE-E | .994 | .032 | .002 (*.000*) | .002 (*.002*) | **.015** (*.014*) |

*Note.* $N = 277$. **Bolded** value exceeds the recommended threshold for equivalence.

*Construct validity*

The CR of the two latent factors (*comprehension and execution* and *emotional regulation*) were found in the MASE-Q to be .87 and .88, and in the MASE-E to be .91 and .84, respectively. Further, the AVE of these latent factors were found in the MASE-Q to be .64 and .70, and in the MASE-E to be .72 and .63, respectively. As CR > .70 for all factors and AVE > .50, convergent validity was achieved for both latent constructs on both scales.

For the MASE-Q, the √AVE were .80 and .84 for *comprehension and execution* and *emotional regulation,* respectively. These were greater than the correlation between the latent factors (.72), so discriminant validity was achieved. However, the MASE-E failed to satisfy the Fornell-Larcker criterion for discriminant validity. Specifically, while the √AVE for *comprehension and execution* (.85) was greater than the correlation between the two factors (.81), the √AVE for *emotional regulation* (.80) was not.

To provide further evidence of the convergent and divergent validity of the self-efficacy factors, Table 5 outlines the correlations between the constructs measured at time point 1, student prior achievement, and exam performance, as well as the descriptive statistics for these variables. As expected, the two constructs (*comprehension and execution* and *emotional regulation*) measuring self-efficacy in both assessment scenarios demonstrate significant, positive correlations with positive emotions (i.e., enjoyment and hope) and achievement measures while correlating negatively with negative emotions (i.e., anxiety and hopelessness) and exam stress. As these results follow established relationships from the literature (e.g. Bandura, 2010; Pekrun & Stephens, 2010), we conclude additional evidence for the validity of the latent self-efficacy constructs. For both self-efficacy factors, the means were higher for quizzes than for the exam.

Table 5. Descriptive statistics and correlations of achievement measures and latent variables at time point 1.

| | *M (SD)* | Range (*Potential*) | α | Skewness | 1 | 2 | 3 | 4 |
|---|---|---|---|---|---|---|---|---|
| *Self-efficacy* (Quiz) | | | | | | | | |
| 1. Comprehension and execution | 70.99 (*15.53*) | 1-100 (*0-100*) | .87 | -2.79 | | | | |
| 2. Emotional regulation | 69.56 (*18.95*) | 4.33-100 (*0-100*) | .89 | -4.26 | .62** | | | |
| *Self-efficacy* (Exam) | | | | | | | | |
| 3. Comprehension and execution | 67.49 (*16.49*) | 1-100 (*0-100*) | .91 | -3.77 | .73** | .48** | | |
| 4. Emotional regulation | 65.09 (*18.44*) | 3-100 (*0-100*) | .84 | -2.72 | .57** | .74** | .69** | |
| *Emotions* | | | | | | | | |
| 5. Anxiety (before) | 3.37 (*0.75*) | 1.25-5 (*1-5*) | .67 | -0.30 | -.29** | -.34** | -.44** | -.51** |
| 6. Anxiety (during) | 2.89 (*0.90*) | 1-5 (*1-5*) | .86 | 0.06 | -.30** | -.40** | -.38** | -.47** |
| 7. Enjoyment (before) | 2.93 (*0.74*) | 1-4.75 (*1-5*) | .74 | -1.72 | .35** | .22** | .39** | .41** |
| 8. Enjoyment (during) | 2.97 (*0.73*) | 1-5 (*1-5*) | .73 | -1.65 | .34** | .23** | .37** | .38** |
| 9. Hope (before) | 3.30 (*0.65*) | 1.20-5 (*1-5*) | .78 | -2.64 | .49** | .30** | .50** | .51** |
| 10. Hope (during) | 3.28 (*0.67*) | 1-5 (*1-5*) | .73 | -1.69 | .41** | .28** | .50** | .46** |
| 11. Hopelessness (before) | 2.36 (*0.78*) | 1-4.60 (*1-5*) | .84 | 1.73 | -.40** | -.39** | -.50** | -.54** |
| 12. Hopelessness (during) | 2.59 (*0.80*) | 1-5 (*1-5*) | .83 | 2.32 | -.37** | -.38** | -.47** | -.50** |
| 13. Exam stress | 6.48 (*1.70*) | 2-9 (*1-9*) | - | -4.18 | -.23** | -.25** | -.32** | -.33** |
| *Achievement* | | | | | | | | |
| 14. Prior achievement | 6.72 (*2.12*) | 1-9 (*1-9*) | - | -5.12 | .29** | .19** | .37** | .30** |
| 15. Exam performance | 53.99 (*24.72*) | 0-98.33 (*0-100*) | - | -1.55 | .34** | .20** | .38** | .25** |

*Note.* N = 277. * *p* < .05, ** *p* < .01

*Conclusions (Study 2)*

An initial goal in the development of the MASE was to have a broad range and common set of scale items that could be applied to different assessment scenarios. However, the analysis in Study 2 indicated we were unable to retain all ten items or maintain identical scales in both assessment scenarios. This finding could possibly be interpreted as building on evidence of how extensively context matters when considering self-efficacy (Bandura, 1997), in turn, highlighting the significance of existing measures neglecting to account for specificities in context to sufficiently investigate student beliefs. Furthermore, the differences between the MASE-Q and MASE-E emphasise the necessity of careful consideration when context-specific self-efficacy is being researched.

Study 2 confirmed both models successfully demonstrated equivalence over time in a single cohort of students through invariance at the configural, metric, and scalar levels. We finally looked for evidence of construct validity of the scales using data from the first time point. We found convergent validity of the constructs was supported for both the MASE-Q and MASE-E through examining the AVE and CR. This was additionally supported by the relationships between self-efficacy factors with achievement emotions and academic achievement. We found the latent constructs of the MASE-Q demonstrated discriminant validity, but failed to satisfy the Fornell-Larcker criterion for the MASE-E. However, we note that this criterion is recognised as overly-sensitive (or biased) toward rejecting discriminant validity in cases where it exists (Rönkkö & Cho, 2020). The recommendation at this stage is to assess the discriminant validity again on a different set of data, which we present in Study 3. From this study, we conclude the MASE-Q to be a valid and reliable measure over time, with the discriminant validity of the MASE-E needing further investigation before the same conclusion can be drawn.

**Study 3**

*Purpose*

The aim of this study was to establish the invariance of the MASE-E between cohorts of students. Additionally, following the results of Study 2, we sought to check for discriminant validity of constructs in the MASE-E on a new set of data.

*Participants*

Participants were two cohorts of students from the same course described in Study 2. The first cohort of students consisted of data from the first time point in Study 2, collected in the second semester of 2020. The data collection for the second cohort took place in the same course during the first semester of 2021 (March - June). There were 355 students enrolled in the course at the start of the 2021 semester. Out of these, 329 provided their consent to the use of their data from the course. There were 38 items on the questionnaire, out of these, one participant missed one response and one participant missed ten responses. Missing data were inserted using EM-imputation. The dataset consisted of 184 participants who identified as male, 144 as female, and three as gender diverse. Ethnicity of participants was not collected in this study. An approval to conduct the study from the University of Auckland Human Participants Ethics Committee was granted (approval number UAHPEC21976).

*Procedure*

Data collection for the 2020 cohort is described in Study 2. For students in the 2021

cohort, the questionnaire was distributed online in the first week of the semester. Students were allocated a small amount (0.15%) of their final course grade for questionnaire completion. This study only collected responses to the MASE-E, so we were unable to validate the MASE-Q across cohorts.

*Measures*

*The MASE-E*

The survey instrument for both cohorts included the MASE-E described in Study 2.

*Data analyses*

The same procedures described in the previous studies were employed for data cleaning. Additionally, the same relevant procedures and criteria were applied to determine model fit, to validate the psychometric equivalence of the MASE-E across groups and to determine validity.

*Results*

*Confirmatory factor analysis*

All factor loadings of the 7-item MASE-E for each cohort were good or excellent (range = .73 to .91) – see Appendix D for full statistical details. Table 6 presents confirmatory factor analyses of the two-factor model measuring exam self-efficacy for each of the cohorts. All fit indices suggest a good fit for each cohort and the combined samples. Specifically, the $\chi^2/df$ were "good" (< 3.0), as were the TLI and CFI values ($\geq$ .95), RMSEA (< .06), and SRMR (< .05).

Table 6. MASE-E for two student cohorts.

| Model | $\chi^2$ | df | $\chi^2/df$ | TLI | CFI | RMSEA | (90% CI) | SRMR |
|---|---|---|---|---|---|---|---|---|
| Full Sample ($N$ = 606) | 21.01 | 13 | 1.62 | .996 | .998 | .032 | (.000-.056) | .016 |
| 2020 Cohort ($n$ = 277) | 19.14 | 13 | 1.47 | .992 | .995 | .041 | (.000-.078) | .024 |
| 2021 Cohort ($n$ = 329) | 19.17 | 13 | 1.47 | .995 | .997 | .038 | (.000-.072) | .015 |

*Construct validity*

Continuing from Study 2, we tested construct validity for the MASE-E on a new set of data using the 2021 cohort. The CR of the latent factors (*comprehension and execution* and *emotional regulation*) were found to be .94 and .90, respectively. Further, the AVE of these latent factors were found to be .79 and .75, respectively. As CR > .70 and AVE > .50 for all factors, convergent validity was achieved. The $\sqrt{AVE}$ were .89 and .87 for *comprehension and execution* and *emotional regulation,* respectively. These were greater than the correlation between the latent factors (.82), so discriminant validity was now achieved.

*Measurement invariance between groups*

Following the CFA results demonstrating a good fit for both cohorts, we proceeded to test for measurement invariance between the two groups of students using MGCFA (Table 7). Based on the observed changes in CFI and RMSEA values (range = .000 to .004), the MASE-E demonstrated "strict" invariance (i.e, configural, metric, scalar, and residual); thus allowing for comparisons of latent factor means between cohorts.

Table 7. Measurement invariance of the MASE-E between cohorts.

| Model | CFI Unconstrained | RMSEA Unconstrained | ΔCFI (*ΔRMSEA*) | | |
|---|---|---|---|---|---|
| | | | Metric Invariance | Scalar Invariance | Residual Invariance |
| MASE-E | .996 | .028 | .000 (*.002*) | .001 (*.004*) | .002 (*.000*) |

*Note.* N = 606 (2020 *n* = 277, 2021 *n* = 329).

### *Conclusions (Study 3)*

This final study aimed to determine whether the MASE-E was invariant between two different cohorts of students. Confirmatory factor analysis indicated the model for exam-related self-efficacy fit both cohorts. Using the new set of data we were able to conclude the construct validity of the MASE-E. Invariance testing determined the scale obtained residual invariance between the groups of students, suggesting it is a reliable measure.

### Discussion

Based on three studies, we have formulated the MASE-Q and MASE-E. Previous research has largely neglected to measure context-specific self-efficacy in mathematics, limiting the questions that can be investigated and omitting a potentially important factor in explaining, for example, student performance. The MASE-Q and MASE-E aim to distinguish assessment-related self-efficacy from other forms of academic self-efficacy. This could enable researchers to investigate self-efficacy more accurately in the educational context of assessment.

In developing the scales for Study 1, we prioritised the preservation of items to best capture all the components of assessment-related self-efficacy, which is something that has not been done in previous measures. Though we aimed to develop items that could be applied to different assessment scenarios for a broad range of applications, we were unable to develop a general scale, potentially providing further backing to the claim that self-efficacy is specific to context. For the final scales, we were able to demonstrate both the MASE-Q and MASE-E were invariant over time in Study 2, and invariance of the MASE-E between groups of students in Study 3 (this has yet to be researched for the quiz-related context). Study 2 and 3 also provided evidence of the validity of the scales in measuring self-efficacy. Thus, this research has demonstrated the two measures to be valid and reliable.

### Limitations and future directions

A limitation is that the MASE has only been validated under two assessment scenarios. Though the initial goal was to develop a generalisable scale for assessment self-efficacy, the analysis indicated that context-specificity extends to the scale items used in each scenario. Future research should seek to confirm the scales for similar assessment scenarios. Given the inconsistencies from Study 1 to Study 2, further testing with all 10 items may be indicated in future research under assessment scenarios unique to the ones outlined here.

Another limitation is that the scales were validated in the context of postsecondary mathematics. Due to the purposive sampling nature of Study 2 and Study 3, the generalisability of the scales needs to be tested further. In particular, future studies should investigate the reliability and validity of the MASE in other domains

(beyond mathematics) and using different sampling methods. We predict, however, that with appropriate instructions these scales would be effective for research in other disciplines and at different educational levels (such as secondary school). With future use, the reliability and validity of the scale should be checked. Ultimately, it is important to appropriately account for context when researching self-efficacy. Thus, our recommendation would be for researchers to design their own assessment scenario and combine the MASE with task-specific self-efficacy scales where appropriate to ensure that both context and content are being taken into account to gain a more holistic understanding of students' self-efficacy.

The MASE development is the first part of a broader research project that aims to use the scale in investigating how assessment-related self-efficacy changes over a semester and whether changes in student self-efficacy for short online quizzes are associated with or can influence changes in exam-related self-efficacy. Viewed this way, these scales have numerous applications, including investigating whether interventions aimed to target exam-related self-efficacy are effective, how self-efficacy interacts with other aspects of test-taking such as achievement emotions and performance, and measuring variations between students' self-efficacy in different types of assessment.

Assessments are a central part of the educational process at any level, and, in particular, the fundamental way of gauging student learning at university. Given the association between self-efficacy and performance, it is important to analyse self-efficacy with respect to assessment. The newly developed MASE allows researchers to accurately measure this significant construct. In order to promote student success, methods targeting assessment self-efficacy must be researched and implemented. The MASE-Q and MASE-E provide opportunities for researchers to examine the effectiveness in fostering self-efficacy when innovating in assessment delivery.

**Disclosure statement:** The authors report there are no competing interests to declare.

**Funding details:** This work was supported by the University of Auckland Faculty of Science Research and Development Grant [3720159].

**Data availability statement:** The data that support the findings of this study are openly available in figshare at https://doi.org/10.17608/k6.auckland.c.6212449